\documentclass[generic,preprint,reqno]{imsart}
\pdfoutput=1
\usepackage[breaklinks]{hyperref}
\usepackage[round, sort&compress]{natbib}
\usepackage{url}
\usepackage{amsmath,amsfonts,amssymb,amsthm}
\usepackage{subcaption}
\usepackage{graphicx}
\usepackage{verbatim}
\usepackage{booktabs}
\usepackage[noend]{algpseudocode}
\usepackage{algorithm, mdframed}
\usepackage[normalem]{ulem}
\usepackage{stmaryrd}
\usepackage{nicefrac}

\setattribute{journal}{name}{}

\definecolor{darkred}{rgb}{0.5,0,0}
\definecolor{darkgreen}{rgb}{0, 0.3,0}
\definecolor{darkblue}{rgb}{0,0,0.6}
\definecolor{LightGray}{rgb}{.6,.6,.6}

\newcommand{\defn}[1]{\textbf{#1}}
\newcommand{\defas}{:=}

\newcommand{\Nats}{\Naturals}
\newcommand{\eqdist}{\stackrel{\mathrm{d}}{=}}
\newcommand{\iid}{\stackrel{\text{iid}}{\sim}}
\newcommand{\ind}{\stackrel{\text{ind}}{\sim}}

\def\EM#1{\ensuremath{#1}}
\def\mbb#1{\EM{\mathbb{#1}}}

\def\Naturals{{\EM{{\mbb{N}}}}}

\newcommand{\G}{\mbb{G}}

\newcommand{\ER}{Erd\H{o}s--R\'{e}nyi}

\theoremstyle{plain}
\newtheorem{proposition}{Proposition}[section]

\newtheorem{theorem}[proposition]{Theorem}

\theoremstyle{definition}
\newtheorem{definition}[proposition]{Definition}

\theoremstyle{remark}

\newtheorem{example}[proposition]{Example}

\newcommand{\bsbeta}{\boldsymbol\beta}
\newcommand{\bseta}{\boldsymbol\eta}
\newcommand{\bstheta}{\boldsymbol\theta}

\newcommand{\WW}{\mathbf{W}}

\begin{document}

\begin{frontmatter}
\title{Priors on exchangeable directed graphs}
\runtitle{Priors on exchangeable directed graphs}
\runauthor{D.\ Cai, N.\ Ackerman, and C.\ Freer}

\begin{aug}
	\author{\fnms{Diana Cai}%
    \ead[label=e1]{dcai@uchicago.edu}%
    \ead[label=u1,url]{http://www.dianacai.com}}
\address{University of Chicago\\ Department of Statistics\\ Chicago, IL 60637, USA\\
          \printead{e1}}
	\author{\fnms{Nathanael Ackerman}%
    \ead[label=e2]{nate@math.harvard.edu}%
	\ead[label=u2,url]{http://www.math.harvard.edu/~nate}}
\address{Harvard University\\ Department of Mathematics\\ Cambridge, MA 02138, USA\\
          \printead{e2}}
		  \ \\
\and
\author{\fnms{Cameron Freer}%
    \ead[label=e3]{cameron.freer@gamalon.com}%
    \ead[label=e4]{freer@borelian.com}%
    \ead[label=u3,url]{http://www.cfreer.org}}
\address{Gamalon Labs\\Cambridge, MA 02142, USA\\
          \printead{e3}\\
		  and\\
	Borelian Corporation\\Cambridge, MA 02139, USA\\
          \printead{e4}}
\end{aug}

\begin{abstract}
Directed graphs occur throughout statistical modeling of networks,
and exchangeability is a natural assumption when the ordering of vertices does
not matter.
There is a deep structural theory for
exchangeable \emph{un\-directed} graphs,
which extends to the directed case via
measurable objects known as \emph{digraphons}.
Using digraphons,
we first show how to construct models for exchangeable directed graphs,
including special cases such as tournaments, linear orderings, directed
acyclic graphs, and partial orderings.
We then show how to construct priors on digraphons via
the \emph{infinite relational digraphon model} (di-IRM),
a new Bayesian nonparametric block model for exchangeable directed graphs,
and demonstrate inference
on synthetic data.
\end{abstract}



\end{frontmatter}

\maketitle

\section{Introduction}

Directed graphs arise in many applications involving pairwise relationships
among objects, such as friendships, communication patterns in social
networks, and logical dependencies \citep{wasserman1994social}.
In machine learning, latent variable models are popular tools for modeling
relational data in applications such as
clustering \citep{MR883333, C.Kemp:2006:53fd9, DBLP:conf/mlg/XuTYYK07, DBLP:journals/jmlr/AiroldiBFX08},
feature modeling \citep{MR1951262, DBLP:conf/nips/MillerGJ09, DBLP:conf/icml/PallaKG12}, and
network dynamics \citep{DBLP:conf/icml/FuSX09, DBLP:conf/nips/BlundellHB12,
DBLP:conf/icml/HeaukulaniG13, DBLP:conf/nips/0002L13}.

\begin{figure}[t]
    \centering
    \begin{subfigure}[b]{\linewidth}
    \centering
    \hspace{25pt}
    \includegraphics[scale=0.340]{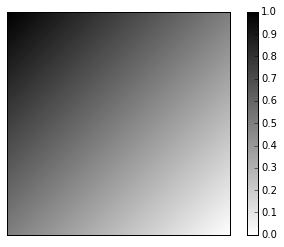}
    \includegraphics[scale=0.340]{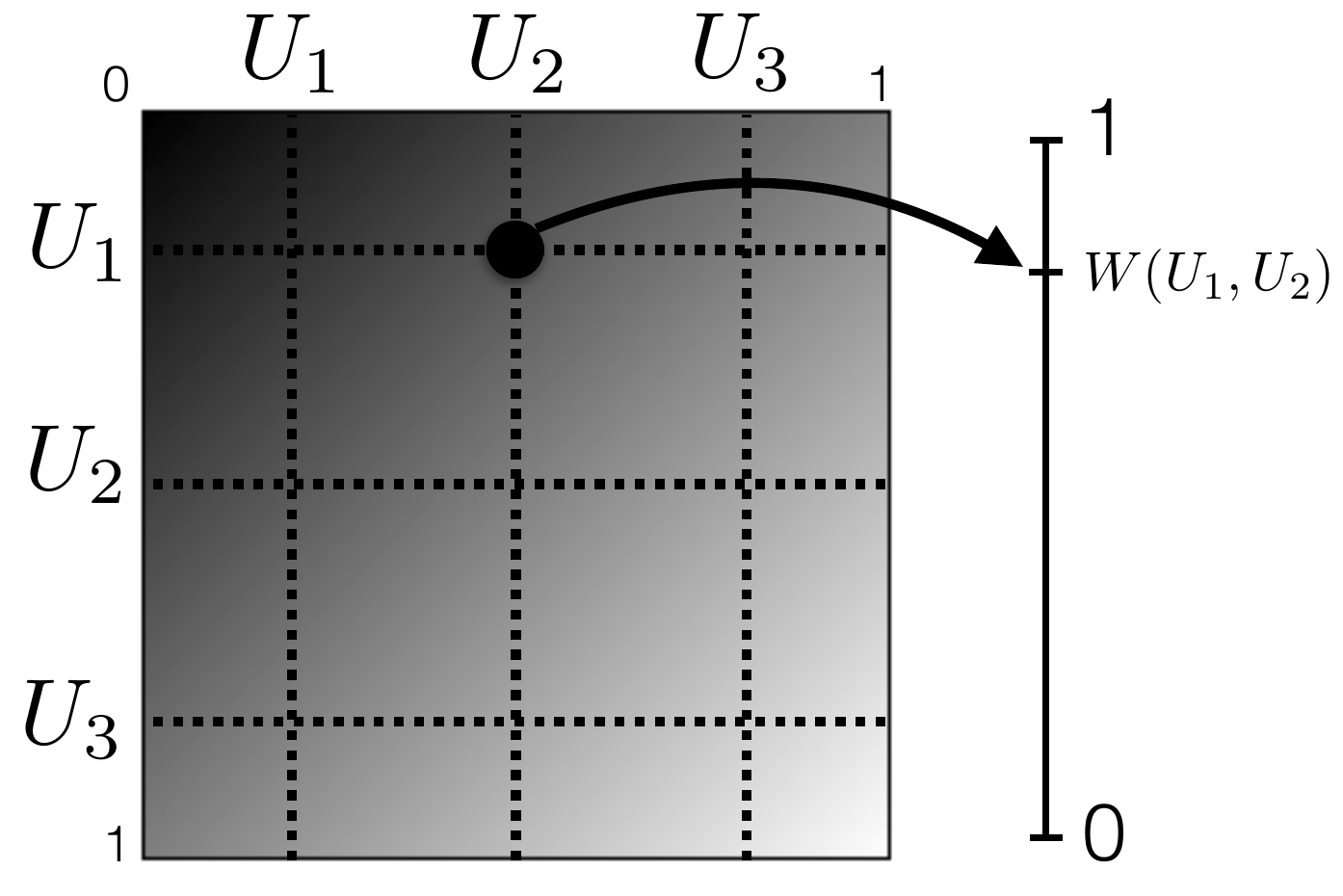}
    \caption{\textbf{left}: Gradient graphon $W$;
        \textbf{right}: Schematic of sampling procedure}
    \end{subfigure}
    \begin{subfigure}[b]{0.23\linewidth}
    \centering
    \includegraphics[scale=0.34]{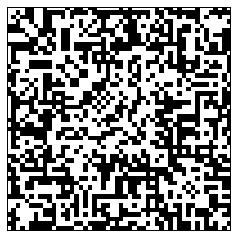}
    \includegraphics[scale=0.34]{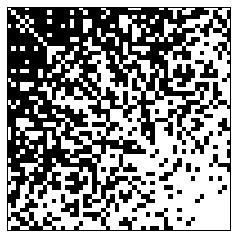}
    \caption{$\G(50,W)$}
    \end{subfigure}
    \begin{subfigure}[b]{0.23\linewidth}
    \centering
    \includegraphics[scale=0.34]{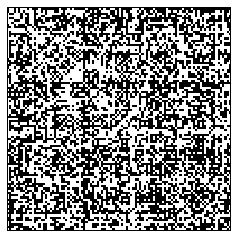}
    \includegraphics[scale=0.34]{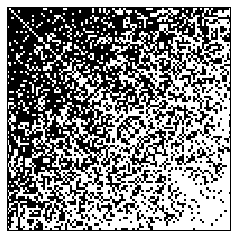}
    \caption{$\G(100,W)$}
    \end{subfigure}
    \begin{subfigure}[b]{0.23\linewidth}
    \centering
    \includegraphics[scale=0.34]{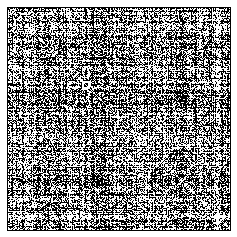}
    \includegraphics[scale=0.34]{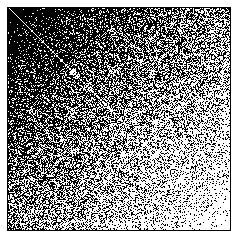}
    \caption{$\G(500,W)$}
    \end{subfigure}
    \caption{(a) An example of a graphon, given by the function
        $W(x,y) = \frac{(1-x)+(1-y)}{2}$.
        (b-d) \textbf{top}: Samples from the finite random graphs $\G(50,W)$, $\G(100,W),$ and $\G(500,W)$,
		shown as ``pixel pictures'' of the adjacency matrix, where black corresponds to $1$ and white to $0$;
        \textbf{bottom}: The samples resorted by increasing order of
        the sampled uniform random variables $U_i$.
        }
    \label{fig-graphon}
\end{figure}

Many such models assume \emph{exchangeability,} i.e., that the
joint distribution of the edges is invariant under permutations of the vertices.
\emph{Undirected} exchangeable graphs have been extensively studied.
The foundational Aldous--Hoover theorem \citep{MR637937, Hoover79}
characterizes
undirected
exchangeable graphs in terms of
certain
measurable functions.
Our perspective in this paper is closer to the equivalent characterization in
terms of \emph{graphons} due to \citet{MR2274085}.
A \defn{graphon} is a symmetric, measurable function $W\colon[0,1]^2 \rightarrow[0,1]$.
Given a graphon $W$, there is an associated countably infinite exchangeable
graph $\G(\Nats,W)$
with random adjacency matrix $(G_{ij})_{i,j\in\Nats}$ defined as follows
(see Figure~\ref{fig-graphon}):
\begin{equation}
\begin{aligned}
    \label{eq-graphon-sample}
	U_i &\iid \text{Uniform}[0,1] \text{~for~} i\in \Nats,\\
    G_{ij} \,|\, U_i, U_j &\ind \text{Bernoulli}(W(U_i, U_j)), \text{~for $i<j$},
\end{aligned}
\end{equation}
and set $G_{ji} = G_{ij}$ for $i<j$, and $G_{ii} = 0$.
Every exchangeable undirected graph can be written as
a mixture of such sampling procedures.
For $n\in\Nats$, we write $\G(n,W)$ to denote the finite random undirected graph on underlying set $\{1, \ldots, n\}$  induced by this sampling procedure.
For more details on graphons and exchangeable graphs,
see the survey by
\citet{MR2463439}
and book by
\citet{MR3012035}.

Most work involving priors on exchangeable graphs has focused on undirected
graphs; for various extensions, see the end of Section~\ref{related-work}.
For directed graphs, much of the work
has extended the undirected case
by using a single
\emph{asymmetric}
measurable function $W_{\text{asym}}\colon[0,1]^2 \rightarrow [0,1]$ to model the
directed graph (see \citet[\S4]{OrbanzRoy15} for a survey of such models).
While such an asymmetric function is appropriate for exchangeable bipartite graphs
\citep{MR2463439},
this representation cannot express all exchangeable directed graph models
(see Section~\ref{ssec-asym}).
Exchangeable \emph{directed} graphs are also characterized by a sampling
procedure given by the Aldous--Hoover theorem.
As with the undirected case, we will work with an equivalent formulation
in terms of
measurable objects known as \emph{digraphons} \citep{MR2463439}; see also
\citet{MR2713257}, \citet{MR3068004}, and \citet{2014arXiv1412.8084A}.
The Aldous--Hoover theorem implies that exchangeable directed graphs
are determined by specifying a distribution on digraphons.
Indeed, a digraphon is a more complicated representation for exchangeable
directed graphs than a single asymmetric measurable function;
a digraphon describes
the possible directed edges
between each pair of vertices
\emph{jointly}, rather than independently.
We define digraphons in Section~\ref{background};
for related work, see Section~\ref{related-work}.

\subsection{Contributions}

This paper presents two main contributions.
We first show how digraphons can be used to model
directed graphs, highlighting special cases that
make use
of dependence
in the edge directions.
In particular, we
characterize the form of digraphons that produce
tournaments, linear orderings, directed
acyclic graphs, and partial orderings (Section~\ref{priors}).
We briefly discuss how these formulations can be used to produce
estimators for directed graph models
(Section~\ref{digraphon:est}).

Next, we given an explicit example of a prior on digraphons:
we present the \emph{infinite relational digraphon model} (di-IRM),
a Bayesian nonparametric block model for exchangeable directed graphs, which
uses a Dirichlet process stick-breaking prior to partition the unit
interval and Dirichlet-distributed weights for each pair of classes in the
partition (Section~\ref{di-IRM}).
We derive a collapsed Gibbs sampling inference procedure
(Section~\ref{sec-inference}), and demonstrate applications of inference on
synthetic data (Section~\ref{experiments}), showing some limitations of using the infinite relational
model with an asymmetric measurable function to model edge
directions independently.

\section{Background}
\label{background}

We begin by defining notation and providing relevant background on directed
exchangeable graphs. Our presentation largely follows \citet{MR2463439}.

\subsection{Notation}
Let $[n] \defas \{1,\ldots,n\}$.
For a directed graph (or \emph{digraph}) $G$ whose vertex set $V$ is $[n]$ or $\Nats$, we write
$(G_{ij})_{i,j\in V}$ for its adjacency matrix, i.e.,
$G_{ij}=1$ if there is an edge from vertex $i$ to vertex $j$, and 0 otherwise. We will omit mention of the set $V$ when it is clear.
In general,
for a directed graph,
$(G_{ij})$ may be asymmetric,
and we allow self-loops,
which correspond to values $G_{ii} = 1$ on the diagonal.
The adjacency matrix of an undirected graph (without self-loops) is a symmetric array $(G_{ij})$
satisfying $G_{ii} = 0$ for all $i$.

We write $X\eqdist Y$  to denote that the random variables  $X$ and $Y$ are equal in distribution.

\subsection{Exchangeability for directed graphs}

A random (infinite) directed graph $G$ on $\Nats$
is \defn{exchangeable} if its joint distribution is invariant under all permutations $\pi$
of the vertices:
\begin{align}
    (G_{ij})_{i,j\in\mathbb{N}} \eqdist (G_{\pi(i)\pi(j)})_{i,j\in\mathbb{N}}.
\end{align}
By the Kolmogorov extension theorem, it is equivalent to ask for this to hold only for those permutations $\pi$ that move a finite number of elements of $\Nats$.

Such an array $(G_{ij})$ is sometimes called \emph{jointly exchangeable}. The case where the distribution is preserved under permutation of each index separately, i.e., where
$(G_{ij})\eqdist (G_{\pi(i)\sigma(j)})$
for arbitrary permutations $\pi$ and $\sigma$, is called \emph{separately exchangeable}, and arises for adjacency matrices of bipartite graphs.

\subsection{Digraphons}
\label{ssec-digraphon}

As described by
\citet{MR2463439},
using the Aldous--Hoover theo\-rem one may show that
every exchangeable countably infinite directed graph is
expressible as a mixture of $\G(\Nats, \WW)$ with respect to some distribution on digraphons $\WW$.

We now define digraphons; in Section~\ref{ssec-digraphon-sample} we will describe the sampling procedure that yields $\G(\Nats, \WW)$.

\begin{definition}
A \defn{digraphon} is a 5-tuple
$\WW \defas (W_{00}, W_{01}, W_{10}, W_{11}, w)$,
where $W_{ab}\colon [0,1]^2 \rightarrow [0,1]$, for $a,b\in\{0,1\}$,
and
$w\colon [0,1] \rightarrow\{0,1\}$ are measurable functions
satisfying the following conditions for all $x,y \in [0,1]$:
\vspace{-10pt}
\begin{align}
    \label{eq-symmetry}
    W_{00}(x,y) &= W_{00}(y,x); \notag\\
    W_{11}(x,y) &= W_{11}(y,x); \\
    W_{01}(x,y) &= W_{10}(y,x); \notag
\end{align}
\vspace{-12pt}
and
\vspace{-10pt}
\begin{align*}
W_{00}(x,y) + W_{01}(x,y) + W_{10}(x,y) + W_{11}(x,y) = 1.
\end{align*}

Given a digraphon $\WW$, write $\WW_4$
for the map $[0,1]^2 \rightarrow [0,1]^4$ given by
$(W_{00}, W_{01}, W_{10}, W_{11})$.
\label{def-digraphon}
\end{definition}

The functions $W_{ab}$ represent the joint probability
of $G_{ij}=a$ and
$G_{ji}=b$ for $a,b \in \{0,1\}$,
i.e.,
\begin{align}
    \label{digraphon-joint}
    \Pr(G_{ij} = a, G_{ji} = b) = W_{ab}(U_i,U_j),
\end{align}
conditioned on $U_i$ and $U_j$.
In this way, $W_{00}$ determines the probability of having neither edge
direction between vertices
$i$ and $j$,
$W_{01}$ of only having a single edge to $j$ from $i$ (``right-to-left"),
$W_{10}$ of a single edge from $i$ to $j$ (``left-to-right"), and
$W_{11}$ of directed edges in both directions between $i$ to $j$.
The function $w$ represents the probability of $G_{ii}$;
because it is $\{0, 1\}$-valued, this merely states whether or not $i$ has a self-loop.

(There is an equivalent alternative set of objects that may be used to specify an exchangeable digraph,
where $W_{00}, W_{01}, W_{10}, W_{11}$
are as before and $p \in [0,1]$ gives
the marginal probability of a self-loop, which is independent of the other edges;
see \citet{MR2463439} for details.)

\subsection{Sampling from a digraphon}
\label{ssec-digraphon-sample}

\begin{figure}
    \centering
    \includegraphics[scale=0.4]{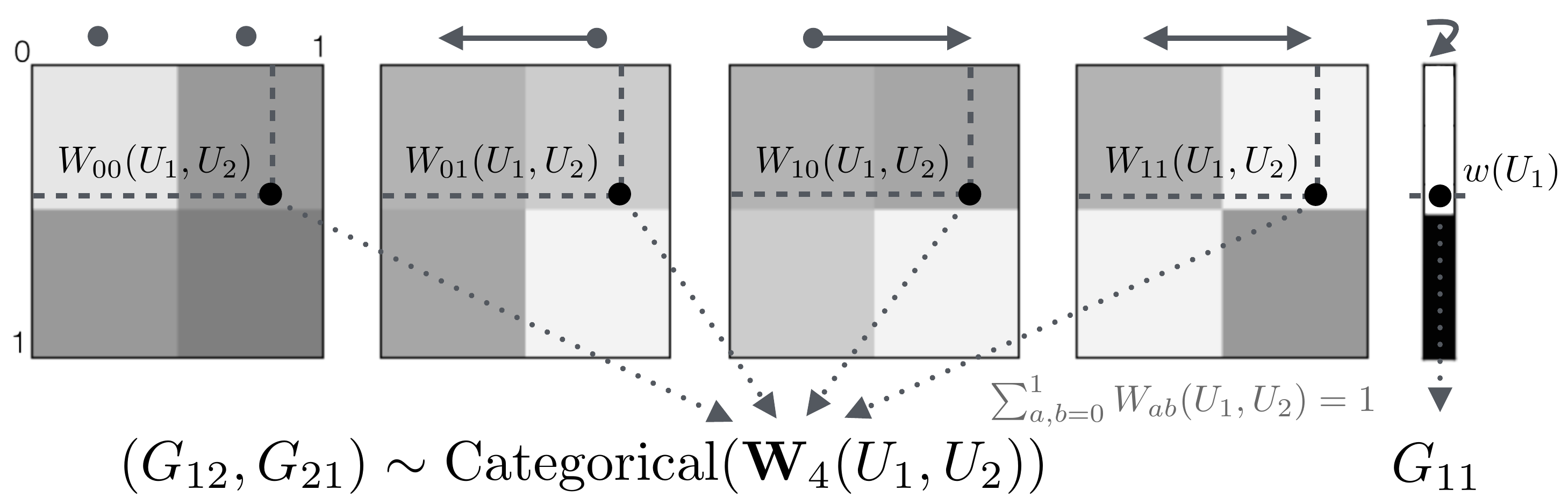}
    \caption{Schematic illustrating digraphon sampling procedure for
         $\WW = (W_{00},W_{01},W_{10},W_{11},w)$.
		 The $x$-axis is vertical and $y$-axis horizontal, with $(0,0)$ in the
         upper left, so that the notation $W_{ab}(x,y)$ coheres with the usual $(\text{row}, \text{column})$ convention for matrix indexing.
    }
    \label{fig:digraphon-schematic}
\end{figure}

The adjacency matrix $(G_{ij})_{i, j\in\Nats}$ of a countably infinite random
graph $\G(\Nats,\WW)$ is determined by the following sampling procedure:
\begin{enumerate}
    \item Draw $U_i \iid \text{Uniform}[0,1]$ for $i\in \Nats$.
    \item For each pair of distinct vertices $i,j$,
            assign the edge values for $G_{ij}$ and
                $G_{ji}$ according to an independent
				$\text{Categorical}(\mathbf{W}_4(U_i, U_j))$
				such that
                Equation~\eqref{digraphon-joint} holds.
        \item Assign self-loops $G_{ii} = w(U_i)$ for all $i$.
\end{enumerate}
In other words, in step 2 we assign
$(G_{ij}, G_{ji})\,|\, U_i, U_j \ind \text{Categorical}(\mathbf{W}_4(U_i, U_j))$, where we interpret the categorical random variable as a distribution over the choices $(0,0), (0,1), (1,0), (1,1)$, in that order.
Note that step 2 is well-defined by the symmetry condition in Equation~\eqref{eq-symmetry}.
Figure~\ref{fig:digraphon-schematic} illustrates
this sampling procedure via a schematic.

An analogous sampling procedure yields \emph{finite} random digraphs:
Given $n\in\Nats$,
in step 1, instead sample only $U_i$ for $i \in [n]$.
Then determine $G_{ij}$ for $i,j\in[n]$ as before.
We write $\G(n, \WW)$ to denote the random digraph thereby induced on $[n]$.

\subsection{Aldous--Hoover theorem for directed graphs}
\label{sec-aldous-hoover}

\citet{MR2463439} derived the following corollary of
the Aldous--Hoover theo\-rem for directed graphs.

\begin{theorem}
[Diaconis--Janson] Every exchangeable random countably infinite directed graph is
obtained as a mixture of $\G(\Nats, \mathbf{W})$; in other words, as
$\G(\Nats, \mathbf{W})$ for some \emph{random} digraphon $\mathbf{W}$.
\end{theorem}

Therefore the problem of specifying the distribution of an infinite exchangeable
digraph may be equivalently viewed as
the problem of specifying a distribution on digraphons.

\section{Digraphons and statistical modeling}
\label{priors}

We first motivate the use of digraphons
instead of asymmetric measurable functions
for modeling exchangeable directed graphs.
We then discuss the representations via digraphons for several random structures
which are special cases of directed graphs.
Finally, we discuss how to estimate digraphons, in the context of both
Bayesian and frequentist estimation.

\subsection{Modeling limitations of asymmetric measurable functions}
\label{ssec-asym}

Asymmetric measurable functions $W_{\text{asym}}\colon [0,1]^2 \rightarrow [0,1]$
characterize exchangeable bipartite graphs by the Aldous--Hoover theorem for
\emph{separately} exchangeable arrays; for details see
\citet[][\S8]{MR2463439}.
These functions can also be used to generate and model \emph{directed} graphs (without self-loops)
by considering the edge directions $G_{ij}$ and $G_{ji}$ independently,
i.e., $\Pr(G_{ij} = 1) = W_{\text{asym}}(U_i, U_j)$ for all $i\neq  j$,
conditioned on $U_i$ and $U_j$,
according to the following sampling procedure:
\begin{align*}
    U_i    & \iid \text{Uniform}[0,1] \text{~for~} i\in \Nats,\\
    G_{ij} \,|\, U_i, U_j & \ind \text{Bernoulli}(W_{\text{asym}}(U_i, U_j)), \text{~for $i\neq j$}
    ,
\end{align*}
and $G_{ii} = 0$ for $i\in \Nats$.
Currently priors on these asymmetric functions are popular in Bayesian modeling
of directed graphs, as we note in Section~\ref{related-work}.

Asymmetric measurable functions are also equivalent to the following
special case of the digraphon representation.
Via the above sampling procedure,
every asymmetric measurable function $W_{\text{asym}}$ yields the same directed graph as the
digraphon
    $\WW = (W_{00},W_{01},W_{10},W_{11},w)$
	given pointwise by
\[
	\WW(x,y) = ((1-p)(1-q),\, (1-p)q, \, p(1-q), \, pq, \, 0),
\]
where $p \defas W_\text{asym}(x,y)$ and $q \defas W_\text{asym}(y,x)$.
In particular, conditioned on $x = U_i$ and $y = U_j$, the marginal probability
$p(1-q) + pq = p$ of an edge from $i$ to $j$ and
$(1-p)q + pq = q$ of an edge from $j$ to $i$ are independent.

On the other hand, many common kinds of digraphs are not obtainable from a
single asymmetric function.
Consider the following two classes:
\begin{enumerate}
	\item Undirected graphs: between any two vertices $i$ and $j$, there are either no
		edges ($G_{ij} = G_{ji} = 0$), or edges in both directions
		($G_{ij} = G_{ji} = 1$).
	\item Tournaments: between any two vertices
        $i$ and $j$, there is exactly one directed edge, i.e., $G_{ij} = 1$ or
	$G_{ji} = 1$ but not both.
\end{enumerate}
For digraphs of either of these two sorts,
the directions are correlated, and hence not obtainable from the above sampling procedure for
an asymmetric measurable function, as this procedure generates
$G_{ij}$ and $G_{ji}$ independently.
This demonstrates how the use of an asymmetric measurable function
is poorly suited for graphs with correlated edge directions.
Though constructing a model for general digraphs using the function $W_\text{asym}$
leads to misspecification, one might hope to perform inference nevertheless; however, as we show in
Section~\ref{half-undir-half-tourn-example}, doing so may fail to discern structure that may be discovered through posterior inference with respect to a prior on digraphons.

In contrast to the use of asymmetric measurable functions, where one considers edge
directions independently, with digraphons one considers the edge directions
between vertex $i$ and vertex $j$ \emph{jointly}, as in Equation~\eqref{digraphon-joint}.
Thus, digraphons give a more general and flexible representation for modeling digraphs.

\subsection{Special cases}

We discuss several special cases of directed graphs and specify the form of
their digraphon representations.

\subsubsection{Undirected graphs}
\label{subsec-undirected}

Undirected graphs can be viewed as directed graphs with no self-loops, where
each pair of distinct vertices either has edges in both directions or in neither.
Hence a digraphon that yields an undirected graph is one having no probability in the single edge directions, i.e., such that
$W_{01}= W_{10} = 0$ (or equivalently, $W_{00} + W_{11} = 1$) and $w=0$.
Such a digraph is therefore determined by merely specifying the graphon $W_{11}$, where
$W_{00}=1-W_{11}$ is implicit.

\begin{figure}
    \includegraphics[scale=0.48]{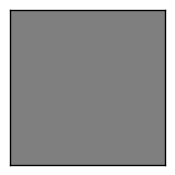}
    \includegraphics[scale=0.48]{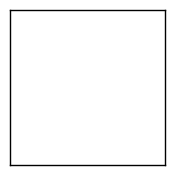}
    \includegraphics[scale=0.48]{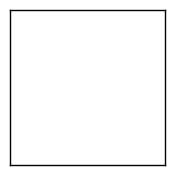}
    \includegraphics[scale=0.48]{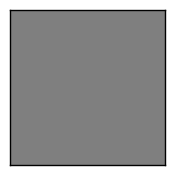}
    \includegraphics[scale=0.48, angle=90]{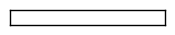}
    \hspace{3pt}
    \includegraphics[scale=0.48]{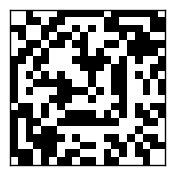}
    \caption{\textbf{left}: \ER\ undirected graph as a digraphon $\WW$;
    \textbf{right}: $\G(20, \WW)$}
    \label{fig-undirected}
\end{figure}

\begin{example}
In Figure~\ref{fig-undirected}, we display an example of a digraphon whose
samples are undirected \ER\ graphs with edge density $\nicefrac12$,
i.e.,
\begin{align*}
	(W_{00}, W_{01}, W_{10}, W_{11}, w) =
    \left(\nicefrac{1}{2}, 0, 0, \nicefrac{1}{2}, 0\right).
\end{align*}
This digraphon corresponds to the graphon $W(x,y) = \nicefrac12$.
\end{example}

\subsubsection{Tournaments}

A \defn{tournament} is a directed graph without self-loops,
where for each pair of vertices, there is an edge in exactly one direction.
In other words, a tournament has $G_{ij} = 1$ if and only if $G_{ji}=0$ for $i \neq j$,
and $G_{ii} = 0$.
Therefore a digraphon yielding a tournament is one satisfying   $w=0$ and $ W_{01} + W_{10} =1$ (or equivalently, $W_{00} = W_{11} = 0$).

\begin{figure}
    \includegraphics[scale=0.48]{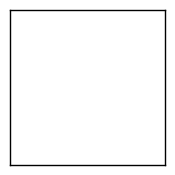}
    \includegraphics[scale=0.48]{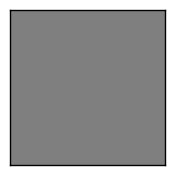}
    \includegraphics[scale=0.48]{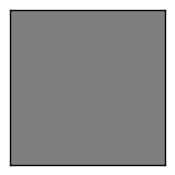}
    \includegraphics[scale=0.48]{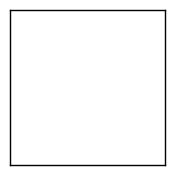}
    \includegraphics[scale=0.48, angle=90]{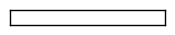}
    \hspace{3pt}
    \includegraphics[scale=0.48]{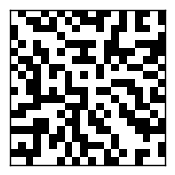}

    \caption{\textbf{left}: Digraphon $\WW$ that yields a generic tournament;
        \textbf{right}: $\G(20, \WW)$}
    \label{fig-tournament}
\end{figure}

\begin{example}
\label{ex-tournament}
An example of a tournament digraphon is displayed in Figure~\ref{fig-tournament}:
\begin{align*}
	(W_{00}, W_{01}, W_{10}, W_{11}, w) =
\left(0, \nicefrac{1}{2}, \nicefrac{1}{2}, 0, 0\right).
\end{align*}
The random tournament induced by sampling from this digraphon is almost surely
isomorphic to a countable structure known as the \emph{generic tournament}.
(For more details on this example, see \citet{MR1106530} and  \citet[Example~9.2]{MR2463439}.)
\end{example}

As discussed in
Section~\ref{subsec-undirected}, exchangeable undirected graphs can be specified
in terms of
single functions $W_{11}$ (graphons) and their associated sampling procedure (described
in Equation~\ref{eq-graphon-sample}).
Similarly,
tournaments also have a
single-function representation and associated
sampling procedure.
Namely, a tournament digraphon is determined by a
measurable function
$W_T: [0,1]^2 \rightarrow [0,1]$
that is anti-symmetric in the sense that $W_T(x,y) = 1 - W_T(y,x)$ for all
$x,y \in [0,1]$
(corresponding to the
digraphon condition
$W_{01}(x,y) = W_{10}(y,x)$).
To sample from $W_T$, first sample $G_{ij}\,|\, U_i, U_j \ind \text{Bernoulli}(W_T(U_i, U_j))$
for $i < j$, and then set $G_{ji} = 1 - G_{ij}$ (and $G_{ii} = 0$).
The digraphon in Example~\ref{ex-tournament}
corresponds to the anti-symmetric, measurable function
$W_T(x,y) = \nicefrac12$.

Tournament digraphons have recently been studied in detail by
\citet{2016arXiv160404271T}, which calls the single function $W_T$ a \emph{tournament kernel}.

Statistical models for tournaments appear in the ranking theory
literature,
often using a variant of the Bradley--Terry model \citep{MR0070925},
first described by \citet{MR1545015}.
For more details, including the relation to graphons, see \citet[\S2.7]{2012arXiv1212.1247C}.
This literature, and related
estimation papers such as \citet{2016arXiv160304556C},
is also often framed in terms of a single-function representation.

\subsubsection{Linearly ordered sets}

A digraph is a (strict) linear ordering when the directed edge relation is
transitive, and
every pair of distinct vertices has an edge in exactly one direction.
Consider the digraphon given by $W_{00} = W_{11} = w = 0$ and $W_{01} = 1- W_{10}$,
where
\begin{align*}
	W_{10}(x,y) =
\begin{cases}
    1 &\text{~if~} x < y, \\
    0 &\text{~otherwise}.
\end{cases}
\end{align*}
The countable directed graph induced by sampling from this digraphon is
almost surely a linear order. In fact, this is essentially the only such example ---
by \citet[\S8]{MR1937832}, its distribution is the same as that of
\emph{every} exchangeable linear ordering.
(In other words, any digraphon yielding the (unique) exchangeable linear ordering is
\emph{weakly isomorphic} to this one; see Section~\ref{experiments} for details.)
Furthermore, the countable linear ordering obtained from sampling this digraphon
is almost surely dense and without endpoints, and hence isomorphic to the
rationals. A finite sample with $n$ vertices has distribution equal to the
uniform measure on all $n!$ ways of linearly ordering $\{1, \ldots, n\}$.

\begin{figure}
        \centering
        \includegraphics[scale=0.48]{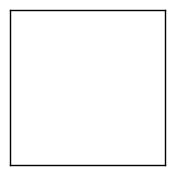}
        \includegraphics[scale=0.48]{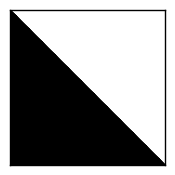}
        \includegraphics[scale=0.48]{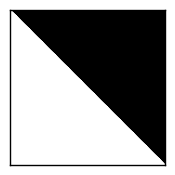}
        \includegraphics[scale=0.48]{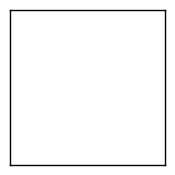}
        \includegraphics[scale=0.48, angle=90]{figs/tourn_w}
        \hspace{3pt}
        \includegraphics[scale=0.48]{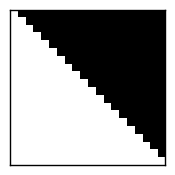}
        \caption{\textbf{left}: Linear ordering digraphon;
        \textbf{right}: $\G(20, \WW)$}
    \label{fig-linear}
\end{figure}

This digraphon is displayed in Figure~\ref{fig-linear} alongside
a 20 vertex random sample,
rearranged by increasing $U_i$;
note that for almost every sample,
the corresponding rearranged graph will have all vertices strictly above the diagonal.

\subsubsection{Directed acyclic graphs}
\label{sec:dag}

A directed acyclic graph (DAG) is a directed graph having no directed path from any vertex to itself.
Various work has focused on models for DAGs (e.g., see \citet{MR2035758}), and especially their use in describing random instances of directed graphical models (also known as Bayesian networks). DAGs also arise naturally as networks describing non-circular dependencies (e.g., among software packages), and in other key data structures.

One can show, using the main result of \citet{MR3324929} (which we describe in Section~\ref{sec:poset}), that
any exchangeable DAG can be obtained from sampling a digraphon satisfying
$W_{10}(x,y) = 0$ for $x \ge y$ and $W_{11} = w = 0$.
Note that this constrains the digraphon to have the same zero-valued regions as those in
the canonical presentation of a linear ordering digraphon (as described
above and displayed in Figure~\ref{fig-linear}), except that
$W_{00}$ may be arbitrary.
(Equivalently,
for $x < y$, the value $W_{10}(x,y)$ may be chosen arbitrarily, so that
the remaining terms are given by $W_{01}(x,y) = W_{10}(y,x)$  and
$W_{00} = 1 - W_{01} - W_{10}$.)
A digraphon of this form thereby specifies one way in which the exchangeable DAG
can be \emph{topologically ordered} (i.e., extended to some exchangeable linear
ordering).

Specifying a digraphon in this way always yields a DAG upon sampling,
as the standard linear ordering on $[0,1]$ does not admit directed cycles,
and one can show
that all exchangeable DAGs arise in this way,
as mentioned above.

\begin{figure}
        \centering
        \includegraphics[scale=0.48]{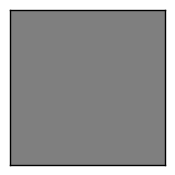}
        \includegraphics[scale=0.48]{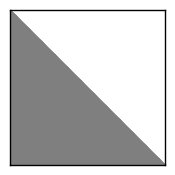}
        \includegraphics[scale=0.48]{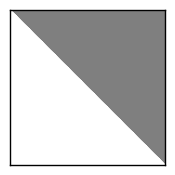}
        \includegraphics[scale=0.48]{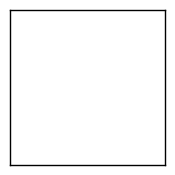}
        \includegraphics[scale=0.48, angle=90]{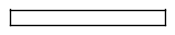}
        \hspace{3pt}
        \includegraphics[scale=0.48]{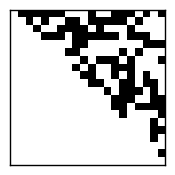}
        \caption{\textbf{left}: Example of a DAG digraphon; \textbf{right}:
            $\G(20, \WW)$}
    \label{fig-dagon}
\end{figure}

\begin{example}
An example of a digraphon that yields exchangeable DAGs
is the generic DAG digraphon given by
\begin{align*}
    W_{00} &= \nicefrac12,\\
    W_{10}(x,y) &=
    \begin{cases}
        \nicefrac12 &\text{~if~} x < y, \\
        0 &\text{~otherwise, and}
    \end{cases}\\
    W_{11} &= 0,
\end{align*}
where $W_{01}$ is such that $W_{01}(x,y) = W_{10}(y,x)$.
This example is displayed in Figure~\ref{fig-dagon}.
We can see that the reordered sample is indeed a DAG, as the edges clearly all
lie above the diagonal in the adjacency matrix.
\end{example}

\subsubsection{Partially ordered sets}
\label{sec:poset}

A partially ordered set, or poset, is a set with a
binary relation $\preceq$ that is
reflexive, antisymmetric, and transitive.
A poset can be viewed as a digraph having a directed
edge from $a$ to $b$ if and only if $a \preceq b$.
Note that the transitive closure of any DAG is a poset, i.e.,
if in a DAG, there is a directed path from $a$
to $b$, the transitive closure has an edge from $a$
to $b$, thereby producing a partial ordering.
(One can similarly define
the ``transitive closure digraphon'' of a digraphon that yields DAGs to obtain
a digraphon yielding the corresponding transitive closures).
Conversely, any poset (with self-loops removed) is already a DAG. Therefore
exchangeable posets are obtainable by some digraphon of the form described in Section~\ref{sec:dag} (except with $w = 1$), though not all such digraphons yield posets. Analogously, representing an exchangeable poset via a digraphon of this form amounts to specifying a \emph{linearization} of the poset.

\citet{MR2886098} develops a theory of poset limits (or posetons) and their relation to exchangeable
posets. By \citet{MR3324929}, any exchangeable poset is given by some digraphon $\WW$ for which $W_{10}(x,y) > 0$ implies that $x < y$, i.e.,
$\WW$
is \emph{compatible} with the standard linear ordering on $[0,1]$.

\begin{figure}
    \includegraphics[scale=0.48]{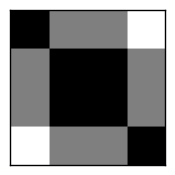}
    \includegraphics[scale=0.48]{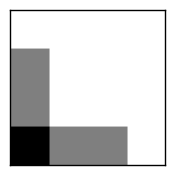}
    \includegraphics[scale=0.48]{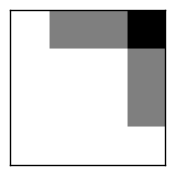}
    \includegraphics[scale=0.48]{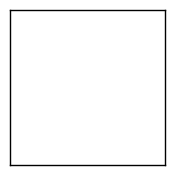}
    \includegraphics[scale=0.48, angle=90]{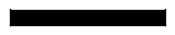}
    \hspace{3pt}
    \includegraphics[scale=0.48]{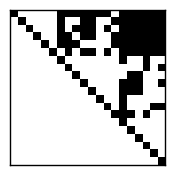}
    \caption{\textbf{left}: Example of a $3\times 3$ SBM poset digraphon; \textbf{right}:
        $\G(20, \WW)$}
    \label{fig-poseton}
\end{figure}

\begin{example}
Consider the following example of a digraphon that yields an exchangeable poset,
specified by the following blockmodel:
\begin{align*}
    W_{10} &=
    \begin{cases}
        \nicefrac12 &\text{~if~} x < \nicefrac14 \text{~and~} \nicefrac14 \leq y  < \nicefrac34 , \\
        \nicefrac12 &\text{~if~} \nicefrac14 \leq x < \nicefrac34 \text{~and~} y \geq \nicefrac34 , \\
		1 &\text{~if~} x < \nicefrac14 \text{~and~} y \geq \nicefrac34 \text{, and} \\
        0 &\text{~otherwise,}
    \end{cases}
\end{align*}
where $W_{11}= 0$,
where $W_{01}$ is such that $W_{01}(x,y) = W_{10}(y,x)$, and
where $W_{00} = 1- W_{01} - W_{10}$.

This example is displayed in Figure~\ref{fig-poseton}.
In particular, the block structure of the model is reflected in the
rearranged sample on the right.
We can see that this is an exchangeable poset:
if the loops (the diagonal) are removed from this digraph, it is
a DAG (as all the edges in the rearranged sample are above the diagonal), and one can check that it is transitively closed.

This is a key example among posets.
Work of \citet{MR0369090} and \citet{MR955579}, characterizing the combinatorial structure of a typical large finite
poset, implies that the sequence of uniform distributions on labeled posets of size $n$ converges
(in the sense of poset limits) to this example.
\end{example}

\subsection{Digraphon estimation}
\label{digraphon:est}

For undirected graphs,  the graphon estimation problem has received considerable
attention in recent years.
In graphon estimation, one seeks to infer either the function $W$, or the
associated
probability matrix with entries $M_{ij} := W(U_i, U_j)$, given a single sample
(or multiple samples) of the graph.
From the Bayesian modeling perspective, one places a prior on
graphons and performs an inference procedure to estimate the parameters of the
random function prior.

From the frequentist perspective, one is interested in producing an estimator for
a fixed graphon, and many such algorithms have been developed, including histogram
and degree-sorting based methods.
To produce a frequentist \emph{digraphon} estimator, one can extend methods
developed for graphons.
Just as a single asymmetric measurable function is insufficient for representing correlated edge directions,
one must likewise estimate the edge directions jointly. Although
a directed graph can be simply represented with a single asymmetric matrix,
digraphon estimators must consider the impact on pairs of entries $(G_{ij},
G_{ji})$ jointly when partitioning, rearranging, or otherwise manipulating
vertices $i, j$.

\paragraph{Histogram estimators for digraphons}
A histogram estimation procedure for graphons partitions the
vertices into several classes, and then uses the average edge density across each pair of
classes as an estimate of the probability of an edge between two vertices in those
classes.
This reduces the problem to that of estimating a partition that yields
a good estimate of these edge densities.
Many methods have been developed for this problem; for further details see the
references within \citet[\S\S1.3 and 1.7]{2015arXiv150806675B}.

To estimate a digraphon (ignoring loops), we must estimate four edge-direction histograms,
where the goal is to estimate a partition of the vertices that
simultaneously yields good estimates of the four types of edge densities.
After producing a partition of the vertices, one likewise computes the
average edge density in each of the four cases, resulting in four histograms.
(If considering loops, there is one additional 1-dimensional histogram whose
estimates are to be jointly optimized by the partition.)

The Frieze--Kannan and Szemer\'edi regularity lemmas lead to bounds on how well
a large graph can be approximated using edge densities across a partition
\cite[Chapters 9 and 10]{MR3012035}; see also \citet{MR1702867}.
The generalization of the Szemer\'edi regularity lemma to directed graphs by
\citet{MR2087940} likewise provides a bound in terms of directed edge densities.

\paragraph{Degree-sorting estimators for digraphons}

Many degree-sorting algorithms have been proposed for graphon estimation.
These algorithms often involve ``sorting'' followed by ``smoothing''.
In the sorting step, the vertices are sorted
by their degrees, where the degree of a vertex $i$ is defined to be $\sum_{j=1}^n G_{ij}$,
In the smoothing step, the $\{0,1\}$-valued adjacency matrix
is used to produce a $[0,1]$-valued matrix using some smoothing algorithm.
For example,
\citet{DBLP:conf/icml/ChanA14} compare a degree-sorting algorithm that uses
total variation distance minimization
as a smoothing step
to one that uses
universal singular value thresholding \citep{2012arXiv1212.1247C} as the
smoothing step.
Degree-sorting estimators
assume that the \emph{degree distribution}
$\int_0^1 W(x,y) dy$ is strictly monotonizable in $x$, i.e., in order for sorting to be
effective, the degrees of the vertices must vary.

This idea can be similarly applied to digraphon estimation:
{First sort} the degrees of the vertices by the four types of edge
directions, to obtain four adjacency matrices, and then smooth these matrices.
It would suffice to require,
after possibly applying a single measure-preserving transformation to $[0,1]$,
that the map $x\mapsto \bigl(\int_0^1 W_{00}(x,y) dy, \int_0^1 W_{01}(x,y) dy, \int_0^1 W_{10}(x,y) dy, \int_0^1 W_{11}(x,y) dy\bigr)$
is strictly increasing with respect to the lexicographic ordering of $[0,1]^4$.

In this paper, we do not comment further on digraphon estimators, but many
other graphon estimation techniques should generalize similarly.
One way of
describing the general pattern is to jointly consider the corresponding
techniques applied to the four matrices obtained from the adjacency matrix restricted to each joint edge type.

\paragraph{Priors on digraphons}
Bayesian approaches may also be use to estimate a digraphon; this is the focus of much of the rest of the paper.
One may likewise use similar techniques to those that have been developed for graphons.
Analogously to
the case of undirected graphs, a Bayesian model for exchangeable directed
graphs involves placing a prior on digraphons.
This is justified by the characterization
in Section~\ref{sec-aldous-hoover}
of exchangeable directed graphs in terms of random digraphons.
We discuss such an approach in depth in Section~\ref{di-IRM}, where we present a Bayesian
nonparametric model based on random partitions using the Dirichlet process.

\section{Infinite relational digraphon model}
\label{di-IRM}

We now proceed to describe a prior on digraphons that makes use of \emph{block structure}.
For directed graphs, the infinite
relational model (IRM) \citep{C.Kemp:2006:53fd9}
models edges between vertices using an asymmetric measurable function and is a
nonparametric extension of the (asymmetric) stochastic block model.
In this section, we present the \emph{infinite relational digraphon model} (di-IRM),
which gives a prior on digraphons.
This model can be viewed as a generalization of the symmetric IRM, a graphon
model, to the digraphon case.
We then show how the di-IRM can be used to model a variety of digraphs, including
ones that cannot be modeled using an asymmetric IRM.

\subsection{Generative model}

We present two equivalent representations of the di-IRM model: (1) a digraphon
representation and (2) a clustering representation.
The digraphon representation uses a stick-breaking Dirichlet process prior to partition the unit
interval, while the clustering representation uses a Chinese restaurant
process prior to partition the vertices.
The difference between the two representations is analogous to that between
the representations of the IRM given by \citet[\S4.1]{OrbanzRoy15}.

\subsubsection{Digraphon representation}
\label{sssection-digraphon}

We first introduce some notation.
Let $\alpha > 0$ be a concentration hyperparameter,
and $\bsbeta \defas (\beta^{(00)}, \beta^{(01)},
\beta^{(10)}, \beta^{(11)})$
be a hyperparameter vector for the weight matrices
$\bseta \defas (\bseta^{(00)}, \bseta^{(01)},
\bseta^{(10)}, \bseta^{(11)})$,
where $\beta^{(ab)} \in [0,\infty)$ for $a,b \in \{0,1\}$.
We allow some (but not all) of the Dirichlet parameters
to take the value zero, at which the corresponding components must be degenerate.
As a shorthand, we write
$\bseta_{r,s} \defas
(\eta^{(00)}_{r,s}, \eta^{(01)}_{r,s},\eta^{(10)}_{r,s}, \eta^{(11)}_{r,s})$
for the 4-tuple of
weights of the classes $r$ and $s$,
where $r,s \in \Nats$.
The following generative process gives a prior on digraphons:
\begin{enumerate}
    \item Draw a partition of $[0,1]$:
        \[ \boldsymbol\Pi \,|\, \alpha \sim \text{DP-Stick}(\alpha).\]
    \item Draw weights for each pair of classes $(r,s)$ of the partition:
    \begin{enumerate}
        \item Draw weights for the upper diagonal blocks, where $r < s$:
        \[
            \bseta_{r,s} \,|\, \bsbeta \sim \text{Dirichlet}(\bsbeta).
        \]
    \item
        Draw weights for the diagonal blocks:
            \begin{align*}
                (\eta_{r,r}^{(00)},
                \eta_{r,r}^{(01)} + \eta_{r,r}^{(10)},
                \eta_{r,r}^{(11)})
				\ |\  \bsbeta &\sim \text{Dirichlet}(\beta^{(00)},\,
				\beta^{(01)}+\beta^{(10)}, \, \beta^{(11)}),
            \end{align*}
                subject to the constraint
            \begin{align*}
				\eta_{r,r}^{(01)} &= \eta_{r,r}^{(10)}.
            \end{align*}
        \item Set weights for the lower diagonal blocks, where $r > s$, such that the symmetry
            requirements in Equation~\eqref{eq-symmetry} are satisfied:
            \begin{align*}
                \eta^{(00)}_{r,s} &= \eta^{(00)}_{s,r}, \quad
                \eta^{(11)}_{r,s} =  \eta^{(11)}_{s,r}, \\
                \eta^{(01)}_{r,s} &= \eta^{(10)}_{s,r}, \quad
                \eta^{(10)}_{r,s} =  \eta^{(01)}_{s,r}.
            \end{align*}
    \end{enumerate}
\end{enumerate}
In Section~\ref{ssec-special} we show different types of random digraphons
that arise from various settings of $\bsbeta$.
The partition is drawn from a Dirichlet stick-breaking process:
for each $i\in\Nats$, draw
$X_i \iid \text{Beta}(1,\alpha)$,
and for every $k \in \Nats$,
set $V_k = X_k \prod_{i=1}^{k-1} (1-X_i)$, so that $\sum_{k=1}^\infty V_k = 1$,
thereby determining a random partition of $[0,1]$.

The self-loops can be specified using the same partition of $[0,1]$,
either with a deterministic $\{0,1\}$-valued function $w$ or a single
weight $p$, as described in Section~\ref{ssec-digraphon}.
For our purposes, we assume $w=0$.
This generative process fully specifies a random digraphon $\WW$, from which random digraphs
$\G(n, \WW)$ can then be
sampled according to the process given in Section~\ref{ssec-digraphon-sample}.

\subsubsection{Clustering representation}
An alternative representation of the generative process for a partition described above can be
formulated directly in terms of clustering: in this generative process,
each vertex $i$ has a cluster assignment $z_i$. This yields an equivalent assignment to that
given by the digraphon formulation if, after sampling the uniform
random variable $U_i$, we assign vertex $i$ to the cluster corresponding
to the class of the partition of $[0,1]$ that $U_i$ belongs to.

Thus, in place of the first step of the generative process given in
the digraphon representation
(Section~\ref{sssection-digraphon}), we draw a partition of the vertices from a
Chinese restaurant process (CRP) (as described in, e.g., \citet{MR883646}):
    $\mathbf z \sim \text{CRP($\alpha$)},$
where each $z_i$ gives the cluster assignment of vertex $i$, and $\alpha >0$ is
a hyperparameter.
The weights $\bseta$ are drawn in the same manner as in the second step of the digraphon
representation of the di-IRM.
Finally, edges are drawn analogously to the general
digraphon
sampling procedure:
$(G_{ij}, G_{ji}) \ind \text{Categorical}(\bseta_{ z_i, z_j})$,
so that Equation~\eqref{digraphon-joint} holds,
where again the Categorical distribution is over the choices
$(0,0), (0,1), (1,0), (1,1)$.

This representation is particularly convenient for performing inference,
especially when using a collapsed Gibbs sampling procedure, as we show in
Section~\ref{sec-inference}.

\subsection{Special cases obtained from the di-IRM}
\label{ssec-special}

\begin{figure}[p]
    \begin{subfigure}{0.65\linewidth}
        \centering
        \includegraphics[scale=0.32]{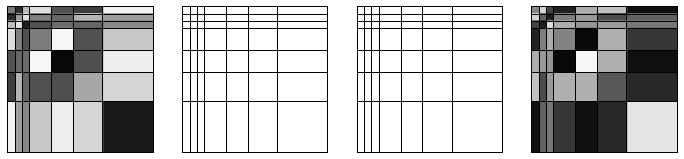}
        \caption{Undirected:
        $\beta^{(01)} = \beta^{(10)} = 0$}
        \label{fig-di-IRM-undir}
    \end{subfigure}
    \begin{subfigure}{0.32\linewidth}
        \centering
        \includegraphics[scale=0.21]{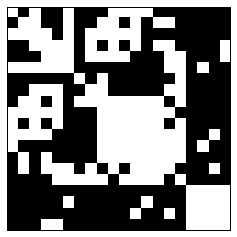}
        \includegraphics[scale=0.21]{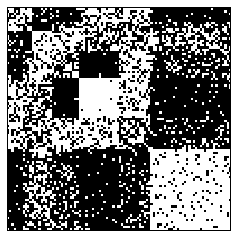}
        \caption{$\G(20,\WW)$, $\G(100,\WW)$}
        \label{fig-di-IRM0}
    \end{subfigure}
    \begin{subfigure}{0.65\linewidth}
        \centering
        \includegraphics[scale=0.32]{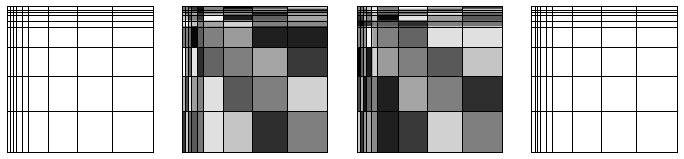}
        \caption{Tournament: $\bsbeta=(0,2,1,0)$}
        \label{fig-dir-IRM-tourn}
    \end{subfigure}
    \begin{subfigure}{0.32\linewidth}
        \centering
        \includegraphics[scale=0.21]{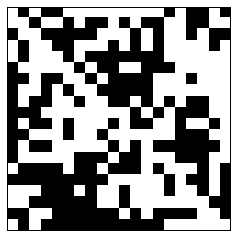}
        \includegraphics[scale=0.21]{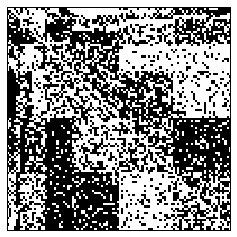}
        \caption{$\G(20,\WW)$, $\G(100,\WW)$}
        \label{fig-di-IRM2}
    \end{subfigure}
    \begin{subfigure}{0.65\linewidth}
        \centering
        \includegraphics[scale=0.32]{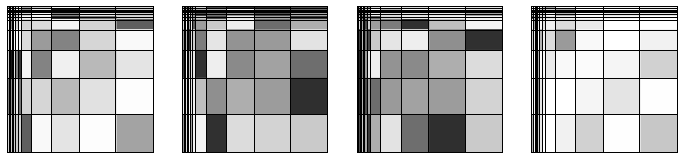}
        \caption{Correlated directions:
        $\bsbeta=(0.9,2,1,0.5)$}
        \label{fig-di-IRM1_fn}
    \end{subfigure}
    \begin{subfigure}{0.32\linewidth}
        \centering
        \includegraphics[scale=0.21]{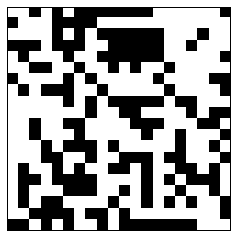}
        \includegraphics[scale=0.21]{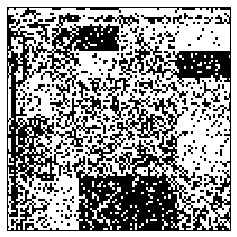}
        \caption{$\G(20,\WW)$, $\G(100,\WW)$}
        \label{fig-di-IRM1}
    \end{subfigure}
    \begin{subfigure}{0.65\linewidth}
        \centering
        \includegraphics[scale=0.32]{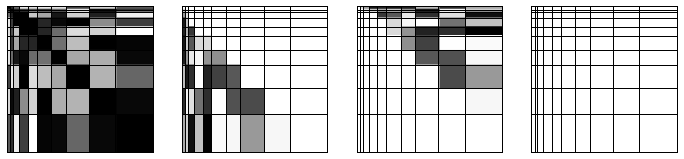}
        \caption{Acyclic:
        $\bsbeta=(0.5,0,0.5,0)$}
        \label{fig-di-IRM_dag_fn}
    \end{subfigure}
    \begin{subfigure}{0.32\linewidth}
        \centering
        \includegraphics[scale=0.21]{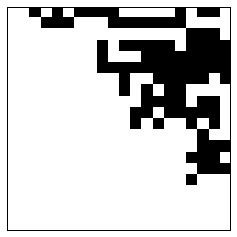}
        \includegraphics[scale=0.21]{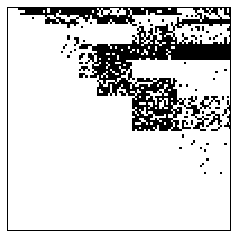}
        \caption{$\G(20,\WW)$, $\G(100,\WW)$}
        \label{fig-di-IRM_dag_fn_samp}
    \end{subfigure}
    \begin{subfigure}{0.65\linewidth}
        \centering
        \includegraphics[scale=0.32]{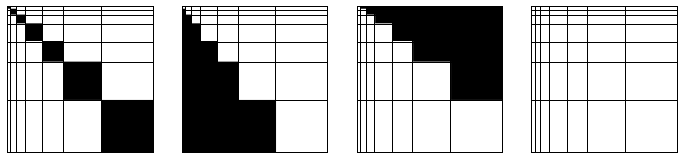}
        \caption{Near-ordering: $\bsbeta=(0,0,1,0)$}
        \label{fig-di-IRM-order}
    \end{subfigure}
    \begin{subfigure}{0.32\linewidth}
        \centering
        \includegraphics[scale=0.21]{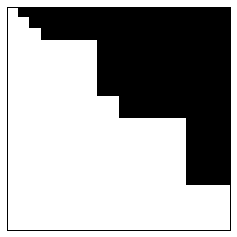}
        \includegraphics[scale=0.21]{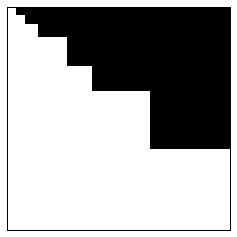}
        \caption{$\G(20,\WW)$, $\G(100,\WW)$}
        \label{fig-di-IRM3}
    \end{subfigure}
    \caption{
        Each row shows a random digraphon drawn from the di-IRM prior
        along with a 20-vertex sample and a 100-vertex sample,
        arranged in order of increasing $U_i$.
		In the smaller samples, one can see certain properties of the digraph
        (e.g., that (b) is symmetric and (d) is a tournament), while in the
        larger samples one can discern block structure with approximate edge densities.
        For (g) and (i), the values given for $\bsbeta$ only apply to
        the classes $r \neq s$, and instead the tuple of hyperparmeters
        $\bsbeta_{r} = (1,0,0)$ is used for the diagonal.
        }
    \label{fig-di-IRM}
\end{figure}

In Figure~\ref{fig-di-IRM}, we display examples of random di-IRM draws
using several settings of the hyperparameter vector $\boldsymbol\beta$.
The parameter settings were specifically chosen to illustrate some of the
special cases the di-IRM model can cover.

\paragraph{Undirected}
To get a prior on graphons using the di-IRM, we can set
$\beta^{(01)} = \beta^{(10)} = 0$.
Figure~\ref{fig-di-IRM-undir} shows a parameter setting that produces undirected
graphs and is equivalent to a symmetric IRM when taking $W_{11}$ to be the IRM;
we can see from the sample on the right that the graph is indeed undirected.

\paragraph{Tournaments}
We can specify a di-IRM tournament prior by setting
$\beta^{(00)} = \beta^{(11)} = 0$.
Figure~\ref{fig-dir-IRM-tourn} shows the parameter setting $\bsbeta = (0, 2, 1, 0)$,
which puts all the mass on the middle two functions.
The tournament structure is easy to see in the 20-vertex sample; for distinct $i$ and $j$, whenever there
is an edge from $i$ to $j$, there is not an edge from $j$ to $i$.

Figure~\ref{fig-di-IRM1_fn} shows a less extreme (non-tournament) variant
that still has strong
correlations between the edge directions,
by virtue of retaining most of the mass on the functions $W_{01}$ and $W_{10}$.
Here we set
$\bsbeta = (0.9, 2, 1, 0.5)$.
Note that the block structure in a sample from this digraphon is more subtle
than in the undirected sample, demonstrating the importance of counting
all four edge-direction combinations rather than just marginals for the two directions.

\paragraph{Directed acyclic graphs}
	To obtain a directed acyclic graph from the di-IRM, we set the hyperparameters so that
	the resulting function $W_{11}$ is empty and $W_{10}$ has nonzero values only on blocks above the diagonal,
	as in Section~\ref{sec:dag}.
	To achieve this,
we set the Dirichlet weight parameters $\bsbeta$ such that
$\beta^{(01)} = \beta^{(11)} = 0$ for the weights $\bseta_{r,s}$ where $r \neq s$,
and for each class $r$
let
$\bsbeta_{r}$
refer to the
3-tuple of hyperparameters
used
for the $\bseta_{r,r}$ weights on the diagonal,
each set to $\bsbeta_r = (1,0,0)$.
With these hyperparameter
choices, we obtain
a directed acyclic di-IRM,
as seen in Figure~\ref{fig-di-IRM_dag_fn}.
We can see in both samples that the directed edges in the resorted sample lie
above the diagonal.
Note that we
make use of $\bsbeta_{r}$ only
in this section, to show how to get a DAG
prior;
in our later inference examples, we use the di-IRM model as introduced in the previous
subsection with the single vector of hyperparameters $\bsbeta$.

\paragraph{Near-ordering}
Consider the hyperparameter settings
	$\bsbeta = (0,0,1,0)$ for the weights $\bseta_{r,s}$
when $r \neq s$,
and $\bsbeta_{r} = (1,0,0)$ for every
class $r$.
The resulting digraph is ``nearly'' ordered,
in the sense that
it
is linearizable and any two elements in different classes are comparable,
as seen in Figure~\ref{fig-di-IRM-order}.
Here $\eta^{(10)}_{r,s}=1$ for any blocks $(r,s)$ above the diagonal, and the
resulting partial ordering is apparent in
both of the resorted samples, with all directed edges above the diagonal.

\subsection{Other partitions for the di-IRM}
Any block model digraphon can be specified in a similar manner: first define a
partition of $[0,1]$, which then gives a partition of $[0,1]^2$; next
let each block on $[0,1]^2$ be piecewise constant such that the symmetry
requirements in Equation~\eqref{eq-symmetry}  are satisfied.

In the case where the number of classes and the size of the classes are fixed
parameters, the directed IRM behaves similarly to some random directed SBM.
In addition to the CRP, we can also consider other partitioning schemes.
Alternatively, one can consider other random partitions of $[0,1]$ as well.
For instance, if one is
interested in power law scaling in the number of clusters (and the sizes of
particular clusters), the Pitman--Yor process \citep{MR1434129} provides a
suitable generalization of the Dirichlet process. It has both a stick-breaking and
urn representation analogous to those for the Dirichlet process.

\section{Related work}
\label{related-work}

The stochastic block model (see \citet{MR718088} and \citet{wasserman1994social})
has been well-studied
in the case of directed graphs
\citep{MR608176, MR883333},
including from a Bayesian perspective
\citep{MR883341, MR2055619, nowicki2001estimation}.
Although working within a restricted class of models, already \citet{MR608176} consider the full joint distribution on edge directions, rather than making independence assumptions.

\begin{figure}
    \centering
    \begin{subfigure}{0.2\linewidth}
        \centering
        \includegraphics[scale=0.5]{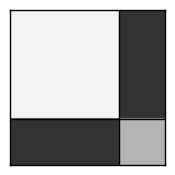}
        \caption{$W_{00}$}
    \end{subfigure}
    \begin{subfigure}{0.2\linewidth}
        \centering
        \includegraphics[scale=0.5]{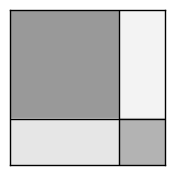}
        \caption{$W_{01}$}
    \end{subfigure}
    \begin{subfigure}{0.2\linewidth}
        \centering
        \includegraphics[scale=0.5]{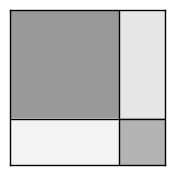}
        \caption{$W_{10}$}
    \end{subfigure}
    \begin{subfigure}{0.2\linewidth}
        \centering
        \includegraphics[scale=0.5]{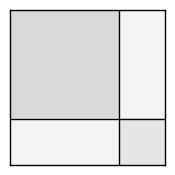}
        \caption{$W_{11}$}
    \end{subfigure}
    \begin{subfigure}{0.1\linewidth}
        \centering
        \includegraphics[scale=0.5, angle=270]{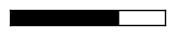}
        \caption{$w$}
    \end{subfigure}
    \caption{Example of a directed stochastic block model digraphon with 2
		classes and a 0.7 division. This example is \emph{assortative}, i.e.,
		there are more edges within the same group than between different
    groups.}
    \label{fig:di-SBM}
\end{figure}

The directed stochastic blockmodel
(di-SBM)
can be represented as a digraphon $\WW_4$ given by four step-functions that are piecewise constant on a finite number of classes.
We display an example of a directed SBM in Figure~\ref{fig:di-SBM}.
The di-IRM model presented in this paper can be seen as a nonparametric
extension of the di-SBM, just as the undirected IRM (introduced independently by
\citet{C.Kemp:2006:53fd9}
and
\citet{DBLP:conf/mlg/XuTYYK07})
is a nonparametric undirected SBM.

Any prior on exchangeable undirected graphs can be described in terms of a corresponding prior on graphons.
As alluded to in the introduction, many Bayesian nonparametric models for graphs admit a nice representation in this form (even if not originally described in these terms); see \citet[\S4]{OrbanzRoy15} for additional details and examples from the machine learning literature, including the IRM.
Likewise, priors on exchangeable digraphs (which have been less thoroughly explored) can be described in terms of the corresponding priors of digraphons, as we have begun to do here.

As noted in \citet{DBLP:conf/nips/LloydOGR12}, when existing models are expressed in these terms, various restrictions (and in particular, unnecessary independence assumptions) become more apparent. As we have seen, the use of the IRM on directed graphs
models the edge directions
as independent (see \citet{2004-019} for examples), a condition that can be straightforwardly relaxed when the model is expressed in the general setting provided by digraphons.

Exchangeable directed graphs have also been considered by \citet{MR2426176}, via an
application of the Aldous--Hoover theorem, although this work does not describe
digraphons explicitly.
We conclude this section by describing several extensions of the graphon
formalism, some of which can be combined with the directed case.
In particular, edges may be more general than $\{0,1\}$-valued.
Variants of graphons for weighted and edge-colored graphs have been considered by
\citet[Chapter~17]{MR3012035} and \citet{MR2426176}.
Graphs with edge multiplicity, or multigraphs,
can be viewed as
integer-valued arrays, a case also covered by the Aldous--Hoover theorem,
although the corresponding extension of graphons is more complicated when the
edge multiplicities are unbounded;
see \citet{MR2754386}, \citet[Chapter~17]{MR3012035}, and \citet{2014arXiv1406.7846K}.
Graphs (that are not necessarily symmetric) with real-valued edges
are also covered by the Aldous--Hoover theorem
through real-valued exchangeable arrays,
and have many applications in statistics and machine learning;
see \citet{DBLP:conf/nips/LloydOGR12} and \citet{OrbanzRoy15}.
The Aldous--Hoover theorem also covers real-valued $d$-dimensional arrays for $d > 2$,
although the corresponding extension of graphons to the case of hypergraphs is
considerably more involved; for details, see
\citet[Chapter~23.3]{MR3012035}, \citet{MR2426176}, and \citet{LOGRhyper}.

\section{Posterior inference}
\label{sec-inference}

In this section, we perform collapsed Gibbs sampling for the di-IRM.
We use the notation for the clustering representation of the di-IRM,
so we can use Gibbs sampling to repeatedly sample the
cluster assignment of each vertex.

Let $G$ be a digraph on $[n]$;
for simplicity we assume that $G$ has no self-edges, and that, as in
Section~\ref{sssection-digraphon}, the di-IRM parameters are chosen so that no
self-edges are produced. Consider a partition of $[n]$ into
a countably
infinite number of clusters, and for
$i\in[n]$, let $z_i\in \Nats$ denote the cluster assignment of $i$.
Write
$\mathbf{z}$ for the vector of all cluster assignments, and
$\bseta$ for the 4-tuple of weight matrices.
Because of
the symmetry requirement of the diagonal,
we
are able to simplify notation as follows:
let
$m_r^* := m_{r,r}^{(01)} + m_{r,r}^{(10)}$, let
$\eta_r^* := \eta_{r,r}^{(01)}+\eta_{r,r}^{(10)}$,
and let
$\beta^* := \beta^{(01)}+\beta^{(10)}$.
Let
$\boldsymbol\beta^* \defas (\beta^{(00)},
\beta^*, \beta^{(11)})$
be the 3-tuple of hyperparameters for the diagonal blocks.

The likelihood of
$G$ being drawn from the di-IRM,
given cluster assignments $\mathbf{z}$ and weights $\bseta$,
is given by
\begin{align*}
p( G \,|\, \mathbf{z}, \boldsymbol\eta) &=
	\prod_{r \leq s}
    \prod_{a,b}(\eta_{r,s}^{(ab)})^{m_{r,s}^{(ab)}}
    \\
        &=
	\Bigl[
	\prod_{r < s}
    \prod_{a,b}(\eta_{r,s}^{(ab)})^{m_{r,s}^{(ab)}}
	\Bigr]
    \,
	\Bigl[
	\prod_{r}
    (\eta_{r,r}^{(00)})^{m_{r,r}^{(00)}}
    (\textstyle\frac12 \eta_{r}^*)^{m_{r}^*}
    (\eta_{r,r}^{(11)})^{m_{r,r}^{(11)}}
	\Bigr]
	,
\end{align*}
where
$m_{r,s}^{(ab)}$ denotes the number of directed edges of type $(ab)$ between
class $r$ and class $s$, for $a, b \in\{0,1\}$ and $r,s\in\Nats$.

Since the weights $\boldsymbol\eta$ have a factorized Dirichlet distribution
prior, we have
\begin{align*}
    p(\bseta \,|\, \bsbeta)&=
	\Bigl[
    \textbf{B}(\bsbeta)^{-1}
    \prod_{r < s}
	\prod_{a,b}
	(\eta_{r,s}^{(ab)})^{\beta^{(ab)}-1}
	\Bigr]
    \\
    & \quad \times
	\Bigl[
    \textbf{B}(\bsbeta^*)^{-1}
    \prod_{r}
	(\eta_{r,r}^{(00)})^{\beta^{(00)}-1}
	(\eta_{r}^*)^{\beta^* -1}
	(\eta_{r,r}^{(11)})^{\beta^{(11)}-1}
	\Bigr],
\end{align*}
where $\textbf{B}(\bstheta) \defas \frac{\prod_i \Gamma(\theta_i)}{\Gamma(\sum_i
\theta_i)}$ is the multivariate beta function.

We sample each cluster assignment $z_i$ conditional on all other assignment variables:
\begin{align}
    \label{cond_dist}
	z_i \,|\, \mathbf{z}_{-i} \sim p(z_i \,|\, \mathbf{z}_{-i}, G) \propto p(G | \mathbf{z}) \, p(z_i | \mathbf{z}_{-i}),
\end{align}
where $\mathbf{z}_{-i}$ denotes the vector of all assignments $z_j$ such that $j \neq i$.

To compute the first term in Equation~\eqref{cond_dist}, we can integrate out the
    parameters $\eta_{r,s}^{(ab)}$:
\begin{align*}
	p(G \,|\, \mathbf{z}) &=
	\Bigl[
    \textbf{B}(\bsbeta)^{-1}
	\prod_{r < s}
	\int
	\prod_{a,b}
	(\eta_{r,s}^{(ab)})^{m_{r,s}^{(ab)}
    +\beta^{(ab)}-1} \, d\eta_{r,s}^{(ab)}
	\Bigr]
	\\
    &
	\quad\,\,
    \times
    \biggl[
    \textbf{B}(\bsbeta^*)^{-1}
    \prod_{r}
    2^{-m_{r}^*} \,
    \\
	&
	\qquad \qquad \qquad
    \int
    (\eta_{r,r}^{(00)})^{m_{r,r}^{(00)}+\beta^{(00)}-1}
    (\eta_{r}^*)^{m_{r}^{*}+\beta^{*}-1}
	(\eta_{r,r}^{(11)})^{m_{r,r}^{(11)}+\beta^{(11)}-1}
     \, d\boldsymbol\eta_{r}
     \biggr]
     \\
    &=
	\Bigl[
	\textbf{B}(\bsbeta)^{-1} \,
    \prod_{r < s}
	\textbf{B}(\textbf{m}_{r,s} + \bsbeta)
	\Bigr]
	\Bigl[
	\textbf{B}(\bsbeta^*)^{-1} \,
    \prod_{r}
    2^{-m_{r}^*} \,
    \textbf{B}(\textbf{m}_{r} + \bsbeta^*)
	\Bigr]
    ,
\end{align*}
where
we simplify calculations on the diagonal using the shorthand
$\textbf m_r \defas (m_{r,r}^{(00)},
m^*_r,
m_{r,r}^{(11)})$,
and
$\boldsymbol\eta_{r} \defas (\eta_{r,r}^{(00)}, \eta^*_r, \eta_{r,r}^{(11)})$.

The second term in Equation~\eqref{cond_dist}
comes from the CRP distribution on $\textbf{z}$:
\begin{align*}
	p(z_i = r \,|\, \mathbf{z}_{-i}) =
    \begin{cases}
		\frac{c_r}{i - 1 + \alpha} \ \text{~if~} c_r > 0,\text{~and}\\
    \frac{\alpha}{i - 1 + \alpha} \ \text{~if~} r \text{~is a new cluster},
    \end{cases}
\end{align*}
where $c_r$ denotes the number of elements in cluster $r$,
and $\alpha > 0$ is the concentration hyperparameter.

We can
reconstruct the weights $\bseta$ using their MAP estimate:
\begin{align}
    \eta_{r,s}^{(ab)} &=
        (m_{r,s}^{(ab)} + \beta^{(ab)})
        /
        N_{r,s},
\end{align}
where
$
N_{r,s} :=  \sum_{a',b' \in\{0,1\}} \bigl(m_{r,s}^{(a'\,b')} +
    \beta^{(a'\,b')}\bigr).
$

\begin{figure}[t]
    \vspace{-5pt}
    \centering
    \begin{subfigure}[b]{\linewidth}
        \centering
        \includegraphics[scale=0.4]{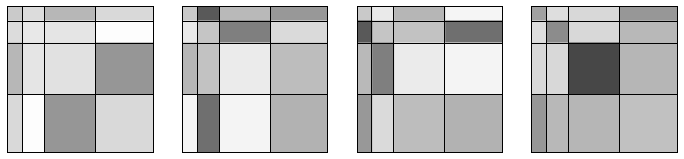}
        \caption{Random digraphon $\WW$ drawn from the di-IRM with $\boldsymbol\beta=(1,1,1,1)$
		}
        \label{subfig:di-IRM-uniform}
    \end{subfigure}
    \begin{subfigure}[b]{\linewidth}
        \centering
        \includegraphics[scale=0.4]{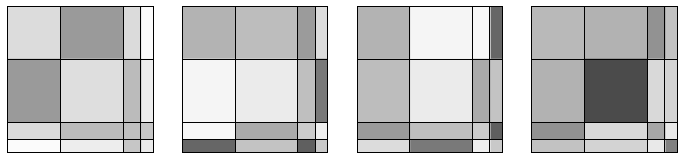}
        \caption{Inferred weights $\bseta$, drawn in proportion to cluster sizes}
    \end{subfigure}
    \begin{subfigure}[b]{0.22\linewidth}
        \centering
        \includegraphics[scale=0.3]{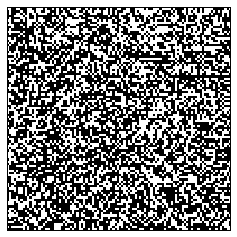}
        \caption{Sampled order}
        \label{subfig:di-IRM-uniform-samp}
    \end{subfigure}
    \begin{subfigure}[b]{0.22\linewidth}
        \centering
        \includegraphics[scale=0.3]{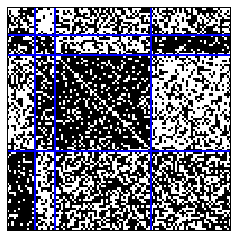}
        \caption{True clusters}
    \end{subfigure}
    \begin{subfigure}[b]{0.22\linewidth}
        \centering
        \includegraphics[scale=0.3]{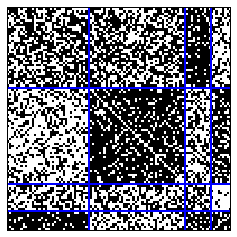}
        \caption{Inferred clusters}
    \end{subfigure}
    \vspace{-8pt}
    \caption{(a) A random digraphon $\WW$ sampled from the di-IRM;
        (b) the inferred weights $\bseta$;
        (c)~a digraph sampled from $\G(100,\WW)$;
        (d) the sample sorted by increasing $U_i$, with true clusters shown via blue lines;
        (e) the sample sorted by the inferred clusters.
        }
    \label{fig:di-IRM-uniform}
    \vspace{-15pt}
\end{figure}

\section{Experiments}
\label{experiments}

In this section, we experimentally evaluate the di-IRM model on synthetic data.
We present two examples: the first is meant to illustrate the correct behavior of
inference on di-IRM parameters, and the second is designed to show the advantage of
using a digraphon representation (given by the di-IRM) over using an asymmetric function (given by the IRM).

Multiple digraphons may induce the same distribution on exchangeable digraphs, in which case they are said to be \emph{weakly isomorphic}. This is not just because a digraphon can be perturbed on a measure-zero set without changing the induced distribution on digraphs, but also because measurable rearrangements of the digraphon will also leave the distribution invariant (analogously to how relabeling the vertices of a digraph does not change it up to isomorphism).
Hence a digraphon $\WW$ is not identifiable from the random digraph $\G(\Nats,\WW)$; in general only its weak isomorphism class can be determined.
For details (in the analogous setting of graphons), see \citet[\S7]{MR2463439} and \citet[\S3.4]{OrbanzRoy15}.

Therefore, in the following inference problems, we can only expect to estimate a digraphon up to its weak isomorphism class.
In a block model, this results in the nonidentifiability of the order of the blocks.

\subsection{Random di-IRM from uniform weights}

We first draw a random di-IRM $\WW$ with the weights $\boldsymbol\beta=(1,1,1,1)$, which is
displayed in Figure~\ref{subfig:di-IRM-uniform}.
We then generate a 100-vertex sample from this digraphon
(Figure~\ref{subfig:di-IRM-uniform-samp}).
We ran a collapsed Gibbs sampling procedure for 200 iterations, beginning from a random initial clustering.
This inference procedure is able to recover the original weights, up to
reordering;
the inferred weight matrices are displayed in
Figure~\ref{fig:di-IRM-uniform}, drawn in proportion to the inferred cluster
sizes.

\begin{figure}[t]
    \centering
    \begin{subfigure}[b]{\linewidth}
        \centering
        \includegraphics[scale=0.5]{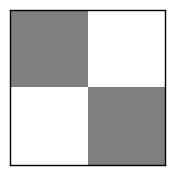}
        \includegraphics[scale=0.5]{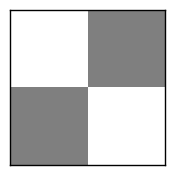}
        \includegraphics[scale=0.5]{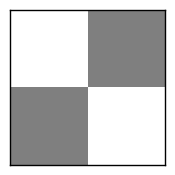}
        \includegraphics[scale=0.5]{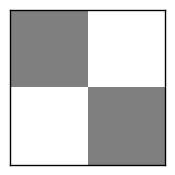}
        \caption{Half-undirected, half-tournament digraphon $W_{00}, W_{01},
        W_{10}, W_{11}$}
        \label{fig:undir_tourn_digraphon}
    \end{subfigure}
    \begin{subfigure}[b]{0.32\linewidth}
        \centering
        \includegraphics[scale=0.78]{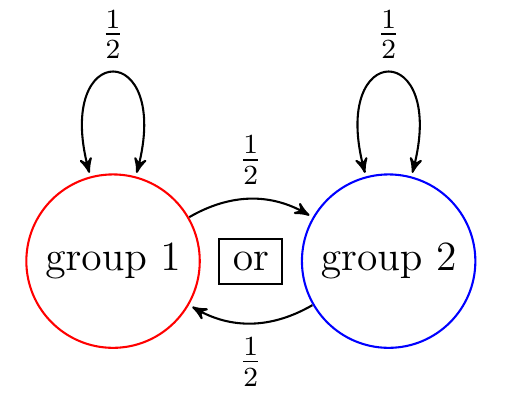}
		\caption{Schematic of example}
        \label{fig:undir_tourn_digraphon_schematic}
    \end{subfigure}
    \begin{subfigure}[b]{0.63\linewidth}
        \centering
        \includegraphics[scale=0.26]{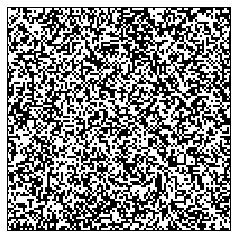}
        \includegraphics[scale=0.26]{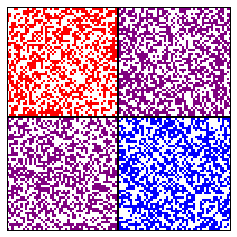}
        \includegraphics[scale=0.26]{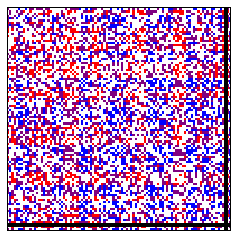}
        \caption{\textbf{left}: Sample; \textbf{middle}: Sample reordered by di-IRM clusters (true
            clusters are colored), \textbf{right}: Sample reordered by IRM clusters}
        \label{fig:undir_tourn_digraphon_sample}
    \end{subfigure}
    \caption{Half-undirected, half-tournament block model example. In the
    schematic, arrows show the probability of connecting in that direction; i.e.,
	any two distinct vertices in the same group have probability
	$\nicefrac12$ of an arrow in both directions (and $\nicefrac12$ of an arrow in neither direction), while for any vertex from group $1$ and vertex from group $2$, either there is
	just an arrow from the first to the second, or there is
	just an arrow from the second to the first, each occurring with probability $\nicefrac12$.
	The bottom right shows the random sample from the
    digraphon and the results of collapsed Gibbs sampling in the di-IRM and the IRM.
    White indicates no edge, red indicates an edge between vertices from
    group 1, blue indicates an edge between vertices from group 2, and
    purple indicates an edge between vertices from different groups.
    Black lines indicate the inferred partition.
}
\end{figure}

\subsection{Half-undirected, half-tournament example}
\label{half-undir-half-tourn-example}

We consider the 2-class step-function digraphon with half the vertices in each class that is given by $w=0$,
\begin{align*}
    W_{00}(x,y) = W_{11}(x,y) =
    \begin{cases}
        \nicefrac{1}{2} &\text{~if~} x < \nicefrac{1}{2} \text{~and~} y<\nicefrac{1}{2}, \\
        \nicefrac{1}{2} &\text{~if~} x \geq \nicefrac{1}{2} \text{~and~} y \geq\nicefrac{1}{2}, \\
        0 & \text{~otherwise,}
    \end{cases}
\end{align*}
and
\begin{align*}
    W_{01}(x,y) = W_{10}(x,y) =
    \begin{cases}
        \nicefrac{1}{2} &\text{~if~} x \geq \nicefrac{1}{2} \text{~and~} y < \nicefrac{1}{2}, \\
        \nicefrac{1}{2} &\text{~if~} x < \nicefrac{1}{2} \text{~and~} y \geq \nicefrac{1}{2}, \\
        0 & \text{~otherwise.}
    \end{cases}
\end{align*}
This digraphon is displayed in Figure~\ref{fig:undir_tourn_digraphon}, and a schematic illustrating the
model is in Figure~\ref{fig:undir_tourn_digraphon_schematic}.
This example demonstrates the importance of being able to distinguish regions having different correlations between edge directions (but the same marginal left-to-right and right-to-left edge probabilities).

We generated a synthetic digraph sampled from $\G(100, \WW)$ and then ran a collapsed
Gibbs sampling procedure for the di-IRM.
We also ran a similar collapsed Gibbs sampler for the IRM.
Both samplers began with a random clustering and ran until the cluster assignments approximately converged.
The results are shown in Figure~\ref{fig:undir_tourn_digraphon_sample}; here the
random sample is displayed alongside the sample resorted according the
clusters inferred using the di-IRM model, as well as the clusters inferred by the
IRM model.
In both resorted images, the true clusters are colored,
white indicates no edge, red indicates an edge between vertices from
group 1, blue indicates an edge between vertices from group 2, and
purple indicates an edge between vertices from different groups.
Note that the true
clusters are correctly inferred using the di-IRM model, as reordering the
vertices according to the inferred clusters identifies the true groups,
while the IRM model fails to discern the correct
structure.
The IRM only considers the
marginal left-to-right and right-to-left edge probabilities, which do not
distinguish the two clusters; in this particular inference run, almost all
vertices were put into the first of the two clusters, which is consistent with
not being able to distinguish between vertices with similar marginal edge
probabilities.
This result is
what one would expect from an algorithm that has inferred uniform independent edge probabilities,
i.e., the edge probabilities of an \ER\ graph.

\section{Discussion}

We have described how priors on digraphons can be used in the statistical modeling of exchangeable dense digraphs, and
have exhibited several key classes of structures that one can model with particular subclasses of these priors.
We have also illustrated why merely using asymmetric measurable functions is insufficient,
as this produces a misspecified model for any exchangeable digraphs having correlations between the edge directions.

While models based on digraphons (and graphons) are almost surely dense (or
empty) and not directly suitable for real-world network applications that are
sparse,
it is still useful to study
models using digraphons (see, e.g., the discussion in \citet[\S7.1]{OrbanzRoy15}).
Some recent work, e.g., \citet{2015arXiv150806675B, 2015arXiv151203099V, 2016arXiv160107134B, nips-Herlau-Mikkel, nips-Cai-Campbell-Broderick, 2016arXiv160304571C, 2014arXiv1401.1137C},
has pointed to methods for extending exchangeable graphs to the case of sparse graphs, but many interesting problems remain.

\section*{Acknowledgments}
A preliminary version of this material was presented at the 10th Conference on Bayesian Nonparametrics in
Raleigh, NC during June 2015, and
an extended abstract was presented at the NIPS 2015 workshop Bayesian Nonparametrics: The Next Generation during December 2015.
The authors would like to thank
Tamara Broderick,
Vikash Mansinghka,
Peter Orbanz,
Daniel Roy, and
Victor Veitch
for helpful conversations,
and Rehana Patel and Daniel Roy for comments on a draft.

This material is based upon work supported by the United States Air Force and the Defense Advanced Research Projects Agency (DARPA) under Contract Numbers FA8750-14-C-0001 and FA8750-14-2-0004.
Work by D.\,C.\ was also supported by a McCormick Fellowship and Bernstein Award at the University of Chicago, an ISBA Junior Researcher Travel Award, and an ISBA@NIPS Special Travel Award.
Work by C.\,F.\ was also supported by Army Research Office grant number W911NF-13-1-0212 and a grant from Google.
Any opinions, findings and conclusions or recommendations expressed in this material are those of the authors and do not necessarily reflect the views of the United States Air Force, Army, or DARPA.

\bibliographystyle{plainnat-mod}
\bibliography{graphons}

\end{document}